\documentstyle[12pt]{article}    % Basic style
\newskip\stdskip                      % standard vertical space
\newenvironment{prf}{\par{\bf Proof}\ }{\hfill$\Box$\par}
\newcommand{\cl}{\centerline}         %  Centerline
\input xy
\xyoption{all}

\newtheorem{thrm}{Theorem}%[section]
\newtheorem{lem}[thrm]{Lemma}
\newtheorem{cor}[thrm]{Corollary}
\newtheorem{rem}[thrm]{Remark}
\newtheorem{defn}[thrm]{Definition}
\newtheorem{exmpl}[thrm]{Example}
\newtheorem{prop}[thrm]{Proposition}

\def\maprt#1{\, \smash{\mathop{\longrightarrow}\limits^{#1}}\, }

\def\bd{  \begin{diagram}    }
\def\ed{  \end{diagram}      }

\def\thm #1  {\medskip\noindent{\bf #1}\quad}
\def\pf #1 {\smallskip\noindent{\bf #1}\par\nobreak\noindent}
\def\term #1{{\bf #1}}

\DeclareSymbolFont{AMSb}{U}{msb}{m}{n}
\DeclareMathSymbol{\ZZ}{\mathbin}{AMSb}{"5A}
\def\s{\Sigma}
\def\smsh{\wedge}

\def\wdg{\vee}

\def\inclds{\hookrightarrow}

\def\of{\circ}
\def\twdl{\widetilde}
\def\<{\langle}
\def\>{\rangle}

\def\sseq{\subseteq}

\def\bprp{\begin{prop}}
\def\eprp{\end{prop}}

\def\bthm{\begin{thrm}}
\def\ethm{\end{thrm}}

\def\blem{\begin{lem}}
\def\elem{\end{lem}}

\def\bcor{\begin{cor}}
\def\ecor{\end{cor}}

\def\brmk{\begin{rem}}
\def\ermk{\end{rem}}

\def\bdfn{\begin{defn}}
\def\edfn{\end{defn}}

\def\bexm{\begin{exmpl}}
\def\eexm{\end{exmpl}}

\def\sp{{\mathrm{Sp}}}
\def\wcat{{\mathrm{wcat}}}
\def\cat{{\mathrm{cat}}}

\def\|{\, \bigm|\, }
\def\sq{{\mathrm Sq}}
\def\d{\overline\Delta}

\def\cat{{\mathrm cat}}

\def\c{{\cal C}}

\def\!{}

\begin{document}
\setlength{\abovedisplayskip}{\stdskip}     % Reduces the space
\setlength{\belowdisplayskip}{\stdskip}     % around displays.
%
%
%                          Title page
%                          ==========
%
%      Acknowledgements should not appear on this page.
%      Please place these at the end of your introduction.
%
%
% This page will be reformatted by the Journal.  It will speed up the
% publication of your article if you type the title page information
% onto the standard form below.
%
% The following macros will turn this information into a basic format for
% the page.  You are welcome to adapt these macros in any way you like.
%
%    Define the various ingredients of the title page:
%
\def\title#1{\def\thetitle{#1}}
\def\authors#1{\def\theauthors{#1}}
\def\author#1{\def\theauthors{#1}}
\def\address#1{\def\theaddress{#1}}
\def\secondaddress#1{\def\thesecondaddress{#1}}

\def\email#1{\def\theemail{#1}}
\def\url#1{\def\theurl{#1}}
\long\def\abstract#1\endabstract{\long\def\theabstract{#1}}
\def\primaryclass#1{\def\theprimaryclass{#1}}
\def\secondaryclass#1{\def\thesecondaryclass{#1}}
\def\keywords#1{\def\thekeywords{#1}}
%
%    Knuth's \ifundefined macro (needed to check for optional items):
\def\ifundefined#1{\expandafter\ifx\csname#1\endcsname\relax}
%
%   Basic title page layout (edit this macro if you
%   wish to adjust the title page layout) :
%
\long\def\maketitlepage{    % start of definition of \maketitlepage

\vglue 0.2truein   % top margin

% title :
%
{\parskip=0pt\leftskip 0pt plus 1fil\def\\{\par\smallskip}{\Large
\bf\thetitle}\par\medskip}

\vglue 0.15truein

% authors :
%
{\parskip=0pt\leftskip 0pt plus 1fil\def\\{\par}{\sc\theauthors}
\par\medskip}%

\vglue 0.1truein

% address(es) email's and URL's (with switches to detect whether the
% optional items have been used) :
%
{\parskip=0pt\small
%{\leftskip 0pt plus 1fil\def\\{\par}{\sl\theaddress}\par}
%\ifundefined{thesecondaddress}\else\cl{and}  % second address?
%{\leftskip 0pt plus 1fil\def\\{\par}{\sl\thesecondaddress}\par}\fi

%\ifundefined{thethirdaddress}\else\cl{and}  % third address?
%{\leftskip 0pt plus 1fil\def\\{\par}{\sl\thethirdaddress}\par}\fi
%\ifundefined{thefourthaddress}\else\cl{and}  % fourth address?
%{\leftskip 0pt plus 1fil\def\\{\par}{\sl\thefourthaddress}\par}\fi

\ifundefined{theemail}\else  % email address?
\vglue 5pt \def\\{\ \ {\rm and}\ \ }
\cl{Email:\ \ \tt\theemail}\fi
\ifundefined{theurl}\else    % URL given?
\vglue 5pt \def\\{\ \ {\rm and}\ \ }
\cl{URL:\ \ \tt\theurl}\fi\par}

\vglue 7pt

{\bf Abstract}

\vglue 5pt

\theabstract

\vglue 7pt

{\bf AMS Classification numbers}\quad Primary:\quad \theprimaryclass

Secondary:\quad \thesecondaryclass

\vglue 5pt

{\bf Keywords:}\quad \thekeywords

%\np  % page break at the end of the title page

}    % end of definition of \maketitlepage
%
%
%   End of macros for basic title page layout
%
%
%  The following lines are for journal use.  Please do not disturb them.
%
%\input gtoutput
%\volumenumber{}\papernumber{}\volumeyear{}
%\pagenumbers{}{}\published{}
%\shorttitle{}  %\shortauthors{}  % for headlines (if needed)
%\proposed{}\seconded{}
%\received{}%\revised{}
%\accepted{}
%
%                  Title page information
%                  ======================
%
%   Type your title page information on the form below following
%   the format of the example.
%
%    \\ is the standard separator (between lines in \title and
%    \address, between authors and email addresses or URL's).
%
%
% Example:  \title{A short spoof paper\\with a two-line title}
% =======   \authors{Albert Einstein\\Leonardo da Vinci}
%           \address{IAS, Princeton}\secondaddress{Renaissance\\Venice}
%           \email{ae@ias.princeton.edu\\ldv@ren.ven.hist}
%           \abstract
%           A short spoof paper with a very short abstract.
%           \endabstract
%           \primaryclass{00-01, 00-02}\secondaryclass{68-00, 68-01}
%           \keywords{Short, spoof, paper}
%
%
%                  Start of title page form
%                  ========================
%
%    Type title, author(s) and address between the curly brackets:

\title{The Lusternik-Schnirelmann\\ Category of $\sp(3)$\\
{\normalsize{\rm Revised Version}}}
\authors{Luc\'ia Fern\'andez-Su\'arez, Antonio G\'omez-Tato,\\
Jeffrey Strom and Daniel Tanr\'e}
\address{}

%  second address, email address and URL (web address), are
%  all optional, uncomment if needed :

\secondaddress{}

%\email{  }
%\url{  }

\abstract   % type your abstract below
We show that the Lusternik-Schnirelmann category
of the symplectic group $\sp(3)$ is $5$.
This L-S category coincides with the cone length
and the stable weak category.
%It follows that $\sp (3)$ satisfies the Ganea
%conjecture -- that is, $\cat(\sp (3)\cross S^n) =
%\cat (\sp(3)) +1$ for all $n\geq 1$.

\endabstract

%  AMS classification numbers, primary and secondary, and keywords :

\primaryclass{55M30}   %LS category
\secondaryclass{22E20} %General Properties and structure of other Lie
%groups

\maketitlepage
%
%%%%%%%%%%%%%%%%%%%%   End of title page
%
%%%%%%%%%%%%%%%%%%%%   Start of main body of article

\section{Introduction}

The Lusternik-Schnirelmann category of a CW complex $X$,
$\cat(X)$, is the
least integer $n$ for which $X$ can be covered by $n+1$ subcomplexes,
each of which is contractible in $X$.
This invariant was introduced by Lusternik and Schnirelmann
in \cite{L-S};  they proved that any smooth function on a compact
manifold $M$ has more than $\cat(M)$ critical points.

Fox showed that Lusternik-Schnirelmann category is a homotopy
invariant \cite{Fox}, so calculating of the category of spaces
is a problem of homotopy theory.   For the analytic
applications, it is crucial to compute the L-S category of %well-known
compact manifolds.  Indeed, determining the L-S category
of Lie groups is the first problem on Ganea's famous problem list \cite{Ganea}.

The present status of this problem is as follows.
Singhof determined the values $\cat({\mathrm{SU}}(n)) = n-1$ and
$\cat({\mathrm{U}}(n))= n$
in \cite{Singhof1}.  The values $\cat({\mathrm{SO}}(n))$ for $n<5$ are easily
computed by classical methods.  More recently, James
and Singhof \cite{J-S} calculated $\cat({\mathrm{SO}}(5)) = 8$.
The symplectic groups have offered the most resistance:  the only known
value is $\cat(\sp(2) ) =3$, which was obtained by Schweitzer \cite{Schw}
%in 1965
using secondary cohomology operations.  Singhof extended
the method to operations of arbitrarily large order, and proved
in \cite{Singhof2} the lower bound $\cat(\sp(n))\geq n+1$.

Actually, Schweitzer proved a little bit more.
He showed that the weak category of $\sp(2)$, which is a
lower bound for $\cat(\sp(2))$, is bounded below by $3$.
Since this lower bound coincides with an upper bound, namely the cone
length, he concluded that $\cat(\sp(2)) = 3$.  (See \S 2 for the definitions.)

Our main result is the following.

\bthm\label{thrm:main}
The symplectic group $\sp(3)$ has Lusternik-Schnirelmann category $5$.
%The Lusternik-Schnirelmann category of the symplectic group $\sp(3)$ is $5$.
\ethm

The proof %of Theorem \ref{thrm:main}
has the same general
outline as Schweitzer's argument for $\sp(2)$.
That is, we show that the cone length
of $\sp(3)$ is at most $5$ and the weak category of $\sp(3)$ is
at least $5$, and conclude that $\cat(\sp(3))=5$.  Our methods are
quite different, however: in place of higher order operations,
we use decompositions of the diagonal map as in
\cite{Strom} together with the Hopf-Ganea invariants studied in \cite{FGT}.

Theorem \ref{thrm:main} gives another proof of the fact,
first proved in \cite{LV},
that $\sp(3)$ satisfies Ganea's conjecture.
To see this, simply apply Theorem 8 in \cite{Strom2}, or observe that,
in the terminology of \cite{Rudyak}, the stable class of
the $5$-fold reduced diagonal $\overline \Delta_5$ is a  detecting class.

The rest of the paper is devoted to the proof of Theorem \ref{thrm:main}.
In \S 2 we establish the basic definitions, notation and results that we
will use.  Section 3 is devoted to an overview of the proof and its reduction to
three propositions.
Their proofs are the content  of the last three sections.

\section{Background}

All spaces are pointed; the basepoint of a space  and the trivial
map between two spaces are both denoted $*$.
The (reduced) cone on a space $A$ will be denoted $\c (A)$.
For a CW complex $X$, we denote by $X_n$ the $n$-skeleton of $X$.
The fat wedge in the $n$-fold cartesian product $X^n$ of $X$ with
itself is the subcomplex
$$
T^n(X) = \{ (x_1, \ldots, x_n) \, |\, {\mathrm at\ least\ one}\
x_i = *\}\sseq X^n.
$$
The cofiber of $T^n(X) \inclds X^n$ is the $n$-fold smash
product of $X$ with itself, denoted $X^{\smsh n}$; the smash
product of two spaces is denoted $X\smsh Y$, as usual.  The diagonal
map $\Delta_n:X\maprt{} X^n$ is the map defined by $\Delta_n(x) = (x,
\ldots , x)$.
The reduced diagonal $\overline \Delta_n: X\maprt{} X^{\smsh n}$
is the composite of $\Delta_n$ with the quotient map $X^n\maprt{} X^{\smsh n}$.

According to Whitehead \cite{Wh}, the \term{category} of $X$ is the least
integer $n$ for which the map $\Delta_{n+1}: X\maprt{} X^{n+1}$ lifts through
$T^{n+1}(X)$, up to homotopy.  The \term{weak category} of $X$,
$\wcat(X)$, is the least $n$ such that $\overline\Delta_{n+1}\simeq *$; clearly
$\wcat(X)\leq \cat(X)$.  See \cite{James} for a survey of
Lusternik-Schnirelmann
category and related invariants.
The \term{strong category}, or \term{cone length}, of $X$ is the least $n$
for which
there is a sequence of cofibrations $A_k\maprt{} X_{k-1}\maprt{}X_k$,
$1\leq k \leq n$,  with $X_0 \simeq *$, $X_n\simeq X$, and $A_k= \s B_k$
for $k>1$ \cite{Cornea,Ganea2}.

We write  $(f,g): X\wdg Y\maprt{} Z$ for the map with components
$f:X\maprt{} Z$ and $g:Y\maprt{}Z$.  We write
$[f,g]:X\maprt{} Y\wdg Z$ for the map with
components $f:X\maprt{}Y$ and $g:X\maprt{}Z$; we use this
notation only in the stable range, so there will be no
ambiguity.

We refer to Toda \cite{Toda} for information about the homotopy
groups of spheres, and we use his notation here.  In particular,
$\omega\in\pi_6(S^3)$ denotes the Blakers-Massey element;
$\omega$ is the attaching map of the $7$-cell in a cellular
decomposition  $S^3\cup D^7\cup D^{10}$ of $\sp(2)$.
Also, for $n>3$, $\nu_n$  is the generator of
the $2$-primary component of $\pi_{n+3}(S^n)$.
Write $C= S^3\cup_\omega D^7$ and
$C_n = S^n\cup_{\nu_n} D^{n+4}$ for $n>4$.

We now recall a result from
\cite{BS} that will be needed in the sequel.

\medskip\noindent{\bf Theorem A}\ \  {\it
If $X$ is a compact $n$-manifold and an H-space, then in any
CW decomposition $X= X_{n-1}\cup_\alpha D^n$,
the map $\alpha$ is stably trivial.}

We will frequently use this fact in the case $X=\sp(n)$, $n=2,3$. It implies that if
the connectivity of $K$ is strictly greater than ${\rm dim}(X)/2$,
then a composite
$X\rightarrow S^{{\rm dim}X}\rightarrow K$
is essential if, and only if, 
$S^{{\rm dim}X}\rightarrow K$
is essential.

The key to our calculation is the following extension of
a result from \cite{FGT}.

\bprp\label{prop:cone}
The space $\sp(3)$ has two cone decompositions:
\begin{enumerate}
\item[{\rm (a)}]
$
\sp(3) = C \cup D^{11} \cup D^{10} \cup D^{14}
\cup D^{18}\cup D^{21},
$
with $\wcat( C\cup D^{11}) = 3$, and
\item[{\rm (b)}]
$
\sp(3) = S^3 \cup \c (C_6) \cup \c (C_{9}) \cup D^{18} \cup D^{21}$.
\end{enumerate}
\eprp

\begin{prf}
Part (a) is proved in \cite{FGT}, as is the fact that
$\sp(3)$ has a cellular decomposition of the form
$$
\sp(3) = S^3 \cup \c (C_6) \cup D^{10}\cup D^{14} \cup D^{18} \cup D^{21}.
$$
It remains to show that the $10$-cell and the $14$-cell can be attached
at the same time as the cone on a map $C_9\maprt{} S^3 \cup \c (C_6)$.
Consider the locally trivial bundle $S^3\maprt{} \sp(3)\maprt{} V_{3,2}$,
where $V_{3,2}$ is the homogeneous space $\sp(3)/S^3$.  This
Stiefel manifold admits a cell decomposition $V_{3,2} =
C_7\cup D^{18}$, so $V_{3,2}$ has
cone length $2$ and hence $\cat(V_{3,2}) = 2$.  We now apply
\cite[Theorem\thinspace 6.2]{Marcum} to complete the proof.
\end{prf}

Notice that it follows from (b) that $\sp(3)$ has cone length
at most $5$.   Therefore Theorem \ref{thrm:main} will be proved
once we show that $\wcat(\sp(3))\geq 5$.

\section{Proof of Theorem \ref{thrm:main}}

Throughout this section let $X$ denote $\sp(3)$ with the CW decomposition
given by Part~(a) of Proposition~\ref{prop:cone}.
Since the subcomplex $A= S^3 \cup \c (C_6)\sseq X$ has
weak category~$2$,
we may decompose the 5-fold reduced diagonal
$\d_5$ using the \emph{main diagram:}
$$
\xymatrix{
X \ar[rr]^(.45){\d_3}\ar@{=}[d] 
   &&  X\smsh X\smsh X \ar[rr]^(.6){1\smsh 1\smsh \d_3}\ar[d]
   &&  X^{\smsh 5}\ar@{=}[d]  
\\
X \ar[rr] \ar@{=}[d] 
   &&  X\smsh X\smsh (X/A) \ar[rr] 
   &&  X^{\smsh 5} 
\\
X \ar[rr] \ar[d]
   && C \smsh C\smsh C_{10} \ar@{^{(}->}[u]\ar[rr]\ar@{=}[d]
   && C^{\smsh 5} \ar@{=}[d]\ar@{^{(}->}[u]
\\
S^{18}\wdg S^{21} \ar[rr]^(.45){q}
\ar@(rd,ld)[rrrr]_{(k,h)}
   && C \smsh C\smsh C_{10} \ar[rr]^(.6){1\smsh 1\smsh
\gamma}
   && C^{\smsh 5} . }
$$
The factorization in the second line is possible because
$\wcat(A) = 2 < 3$, and so $\d_3$
factors through a map $X/A\maprt{} X^{\smsh 3}$.  Neither this map,
nor its restriction $\gamma$ to $C_{10}$, is uniquely determined;
however, the restriction to $S^{10}\sseq C_{10}$ factors uniquely
through the map 
$S^{10}\maprt{} S^{9}$ induced by $\d_3:\sp (2) \maprt{} \sp (2)^{\smsh 3}$,
which is known to be $\eta$ \cite{Schw}.
The third line is a
cellular approximation and the fourth line is obtained by collapsing the 14-skeleton of $X$.

This diagram expresses the reduced diagonal $\d_5$ 
as a composition of several maps, and 
our proof amounts to explicitly determining each of
them.  The details are contained in the proofs of the 
following three propositions.

We first show that the sphere $S^{18}$ does not 
play a role in the map $\d_5$.

\bprp\label{prop:A}
In the main diagram, the composite $S^{18}\maprt{} X^{\smsh 5}$
is trivial.
\eprp 

Thus we are reduced to determining the map 
$S^{21} \maprt{} X^{\smsh 5}$.  Our second result 
shows that this map factors through the 
inclusion $S^{15} \inclds X^{\smsh 5}$. 

\bprp\label{prop:B}
The composite $S^{21} \maprt{} X^{\smsh 5} \maprt{} X^{\smsh 5}/S^{15}$
is trivial.
\eprp

Our final technical result determines the maps
$S^{21}\maprt{}C\smsh C\smsh C_{10}$ and 
$C\smsh C\smsh C_{10}\maprt{} C^{\smsh 5}$.

\bprp\label{prop:C}
There is a homotopy equivalence $$T:C\smsh C\smsh C_{10}
\maprt{} C_{16}\wdg S^{20}\wdg S^{20}$$ such that:
\begin{enumerate}
\item[{\rm (a)}] $T\of q = [ \lambda, 0, 0 ]$ for some map
$\lambda$, and 
\item[{\rm (b)}]  there is a commutative diagram
$$
\xymatrix{
S^{21} \ar[rr]^{\lambda} \ar[rrd]_{\eta} 
   && C_{16} \ar[d] \ar[rr]&&C^{\smsh 5}\ar[d]
\\
&& S^{20} \ar[rr]_(.45){*} && C^{\smsh 5}/S^{15}. }
$$
\end{enumerate}
\eprp 

We now use these three results to prove Theorem \ref{thrm:main}.
The proofs of Propositions \ref{prop:A},
\ref{prop:B} and \ref{prop:C} are given in the 
next three sections.

\smallskip\noindent{\sc Proof  of Theorem \ref{thrm:main}.}\  \
It follows from  Proposition \ref{prop:C}
that the map $h:S^{21}\maprt{} C^{\smsh 5}$ factors through
a map $C_{16}\maprt{} S^{15}$.  An examination of the 
third row of the main diagram reveals that the restriction of 
$C_{16}\maprt{} S^{15}$ to $S^{16}$ is the sixth suspension of
the map 
$S^{10}\maprt{} S^9$ induced by $\d_3$.  This map
has been determined to be $\eta$ in \cite{FGT,Strom,Schw}.
Therefore the map $S^{21}\maprt{} S^{15}$ that lifts $h$
is the composition on the third row of the diagram
$$
\xymatrix{
S^{20} \ar[rr]^{\eta} \ar[d] && S^{19} \ar[d]^{\nu} \\
{*} \ar[d] \ar[rr] && S^{16}\ar[d]\ar[rrd]^{\eta} \\
S^{21} \ar[rr] \ar[rrd]_{\eta}&& C_{16} \ar[d]\ar[rr] && S^{15} \\
&& S^{20}.\\
}
$$
In other words, the map  $S^{21}\maprt{} S^{15}$
is an element of the Toda secondary composition
$\{ \eta, \nu, \eta\}_0$, which is the singleton set
$\{ \nu^2\}$ by \cite[Lemma\thinspace 5.12]{Toda}.

This shows that the 5-fold reduced diagonal 
can be factored as in the diagram
$$
\xymatrix{
X \ar[rr]^{\d_5} \ar[d] && X^{\smsh 5}\\
S^{21}\ar[rr]^{\nu^2} && S^{15}\ar@{^{(}->}[u]_j . }
$$
Since $X$ is an H-space and a manifold, the 
attaching map for the $21$-cell is stably trivial by Theorem A,
and it follows that $j\of \nu^2$ is the {\it unique} map
$S^{21} \maprt{} X^{\smsh 5}$ making the diagram commute. 
Observe that the 22-skeleton of $X^{\smsh 5}$ has the homotopy type of the cofibre
$$
{\bigvee_{i=1}^{5} S^{18}}\maprt{l} S^{15} \maprt{k}
X_{22}= (S^{15}\cup_{2\nu}D^{19})\wdg\bigvee_{i=1}^4
S^{19}\wdg \bigvee_{i=1}^{5} S^{22},
$$
where $l$ is the attaching map $(2\nu, 0,0,0,0)$.
This sequence gives an exact sequence of homotopy groups
$$
\pi_{21}\left({\bigvee_{i=1}^{5} S^{18}}\right)\maprt{l_*}
\pi_{21}\left(S^{15}\right)\maprt{k_*}
\pi_{21}\left( (S^{15}\cup_{2\nu}D^{19})\wdg\bigvee_{i=1}^4
S^{19}\right)\,.
$$
Since $2\nu \of \nu = 0$, it follows that $k\of \nu^2\neq 0$
and so $j\circ \nu^2\neq 0$.
Consequently
$\d_5\not\simeq *$ and hence 
$\cat(\sp(3))=\wcat (\sp(3)) = 5$.
\hfill$\Box$\par

\section{Proof of Proposition \ref{prop:A}}

Since $\wcat (X_{18}) \leq 4$, the map $\d_4:X_{18}\maprt{}X^{\smsh 5}$ 
is trivial.  Also, 
the inclusion $X_{18} \inclds X$ induces the 
inclusion of the first summand $S^{18}\inclds 
S^{18}\wdg S^{21}$ after collapsing $14$-skeleta.  Thus we obtain the commutative
diagrams
$$\xymatrix{
X_{18} \ar[rr] \ar[d] && X \ar[rr] \ar[d] && X^{\smsh 5}\\
S^{18}\ar@{^{(}->}[rr] && S^{18}\wdg S^{21}\ar[rr]_{(k,h)} && 
X^{\smsh 5}\ar@{=}[u] .\\  }\quad
\xymatrix{
{}\ar@{}[d]|
{\mathrm and}\\ {} }\quad
\xymatrix{
X_{18} \ar[rr] \ar[d]  && X^{\smsh 5}\\
S^{18}\ar[rr]_{*}  && 
X^{\smsh 5}.\ar@{=}[u] \\  }
$$
We will show that there is a {\it unique} map
$S^{18}\maprt{} X^{\smsh 5}$ making the second diagram 
commute and conclude that $k=0$.

From the cofiber sequence $X_{18}\maprt{} S^{18}\maprt{}
\s X_{14}\maprt{}
\s X_{18}$, we see that it is enough to show that 
every map $\s X_{14}\maprt{} X^{\smsh 5}$ can be 
extended to $\s X_{18}$.  For this, we consider the 
commutative diagram
$$
\xymatrix{
\s X_{14} \ar[rr] \ar[d] && S^{15} \ar[rr] \ar[d] && X^{\smsh 5}\\
\s X_{18}\ar[rr] && S^{15}\cup_x D^{19}\ar@{-->}[rru] \\  } 
$$
in which the attaching map $x$ is some multiple of $\nu$.
Since the homotopy type of the $19$-skeleton of $X^{\smsh 5}$ is
given by
$$
(X^{\smsh 5})_{19} \simeq 
(S^{15}\cup_{2\nu}D^{19})\wdg\left( \bigvee_{i=1}^4 S^{19}\right), 
$$ 
the extension will be possible if $x$ is an even multiple of 
$\nu$.  If $x$ were an odd multiple of $\nu$, then
the operation $\sq^4: H^{14}(S^{14}\cup_x D^{18})
\maprt{} H^{18}(S^{14}\cup_x D^{18})$ would be an 
isomorphism.  But a simple calculation in 
$H^{*}(\sp(3);\ZZ/2)$ shows that this operation is zero.
\hfill$\Box$\par

\section{Proof of Proposition~\ref{prop:B}}

In this section we write
$\overline X$ for $\sp(3)$ with the CW structure given by Proposition \ref{prop:cone}(b).
We choose a new decomposition of the 5-fold reduced diagonal as follows:
$$
\xymatrix{
\overline X\ar[r]^{\overline\Delta_2} \ar@{=}[d]
   & \overline X \smsh \overline X
\ar[rr]^{\overline\Delta_2\smsh\overline\Delta_3}\ar[d]
   && \overline X^{\smsh 2}\smsh \overline X^{\smsh 3}\ar@{=}[d]
\\
\overline X \ar[r] \ar[d]
   & \overline X/\overline X_3 \smsh \overline X/\overline X_7 \ar[rr]\ar[d]
   && \overline X^{\smsh 2}\smsh \overline X^{\smsh 3}\ar[d]
\\
S^{21}\ar@{=}[d]\ar[r]
    & \overline X/\overline X_7 \smsh \overline X/\overline X_7  \ar[rr]
    &&( {\overline{X}^{\smsh 2}/S^6 )\smsh \overline{X}^{\smsh 3}}
\\
S^{21}\ar[r]^(.29){p}
   & (S^{11}\wdg S^{10}) \smsh (S^{10}\wdg S^{11} )\ar[rr]^{\alpha\smsh \beta}
\ar@{^{(}->}[u]
   && (S^{10}\wdg S^{10})\smsh  S^9.\ar@{^{(}->}[u]
\\ }
$$
The construction of this diagram deserves some explanation.
The second line is obtained from the triviality of
$\overline\Delta_2 |_{\overline X_3}$ and $\overline\Delta_3 |_{\overline X_7}$.
The third line is the result of collapsing the 18-skeleton on the first column and the $7$-skeleton
of the first factor on the last two columns.
The fourth line follows from cellular approximation in each factor.

We determine the maps $p$, $\alpha$
and $\beta$ in the following lemma. 

\blem\label{lem:pab}
The maps $p$, $\alpha$ and $\beta$ are given by
\begin{enumerate}
\item[{\rm (a)}]    $p=[0,\pm 1 ,\pm 1 ]$,
\item[{\rm (b)}]    $\alpha=([\eta,\eta],[\pm 1 ,\pm 1 ])$, and
\item[{\rm (c)}]    $\beta = (\eta, \eta^2)$.
\end{enumerate}
\elem

Assuming these values, we can perform the

\smallskip\noindent{\sc Proof of Proposition \ref{prop:B}.}\  \
Since the manifold $\sp(3)$ is an H-space, the attaching map for the top cell is
stably trivial by Theorem A.  It follows that the map
$\overline{X}\maprt{} (S^{10}\wdg S^{10})\smsh C^{\smsh 3}$, which is in the stable
range, is
trivial if and  only if $(\alpha\smsh\beta)\of p:
S^{21}\maprt{} (S^{10}\wdg S^{10})\smsh S^9$, the extension to $S^{21}$, is
trivial.
When we expand this composite map using the identifications
we have made, we find that:
$$
\begin{array}{rcl}
(\alpha \smsh \beta) \of p
&=&\bigl( ([\eta, \eta],[1,1])\smsh (\eta,\eta^2) \bigr)|_{S^{20}\wdg
S^{21}\wdg
S^{21}}\of [0,1,1]
\\
&=&
( [\eta,\eta] , [\eta^2,\eta^2] , [\eta^2,\eta^2] )\of [0,1,1]
\\
&=&
[ (0 + \eta^2) +  \eta^2, (0 + \eta^2) + \eta^2]= 0.\\
\end{array}
$$

Since $(C^{\smsh 5})_{22}=S^{15}\cup \bigcup_{i=1}^5 D^{19}$, 
the composite $S^{21}\rightarrow C^{\smsh 5}\rightarrow C^{\smsh 5}/S^{15}$
is trivial if each composite
$$
S^{21}\rightarrow C^{\smsh 5}\rightarrow S_i^{19}\vee S_j^{19}
$$
is trivial. The first part of this proof gives the triviality for some $i$ and $j$. We claim
now the triviality for any $i$ and $j$. By Theorem A,
there is a unique map $h$ making  the upper square of the  diagram
$$\xymatrix{
\overline{X}\ar[rr]^{\d_5}\ar[d]&&\overline{X}^{\smsh 5}\\
S^{21}\ar[rr]^{h}\ar[rrd]_{\ast}&&C^{\smsh 5}\ar[d]\ar[u]\\
&&S_i^{19}\vee S_j^{19}
}$$
commutative.
Permutation of the 5 factors of $\overline{X}^{\smsh 5}$ preserves the maps $\d_5$
and $C^{\smsh 5}\rightarrow \overline{X}^{\smsh 5}$. Since $h$ is unique, $h$ is also 
preserved.
The permutation \emph{does} alter the indices $i$ and $j$, and this proves the result.
\hfill$\Box$\par

\smallskip\noindent{\bf Proof  of Lemma \ref{lem:pab}.}\  \
\emph{Proof of {\em (a)}}
By cellular approximation the map $p$ lifts to
$S^{20} \wdg S^{21}\wdg S^{21}$ in the diagram
$$
\xymatrix{
S^{21}\ar[rr]^(.3){p}\ar[rrd]_{[a\eta, \pm 1,\pm1 ]\ \ \, }
    && (S^{11}\wdg S^{10}) \smsh (S^{10}\wdg S^{11} )\\
    && S^{20} \wdg S^{21}\wdg S^{21}.\ar@{^{(}->}[u]\\
}
$$
The
second and third coordinates follow from the cup
product structure in $H^*(\sp(3);\ZZ)$. 
For the first coordinate there are two
possibilities, namely
the Hopf map or the trivial map, that we can represent by $a\eta$
where $\eta\in \ZZ/2$.
We show now that $a=0$.

From the inclusion of the primitive subspace
$\overline{A}=S^3\cup e^7\cup e^{11}\hookrightarrow
\overline{X}$, we construct the commutative diagram
$$\xymatrix{
\overline{X}\ar[rr]^{\d_2}\ar[d]&&\overline{X}\smsh \overline{X}\ar[d]\\
\overline{X}/\overline{A}\ar[rr]\ar[d]&&(\overline{X}/\overline{A})\smsh
(\overline{X}/\overline{A})\\ S^{21}\ar[rr]^-{\overline{H}(\varphi)}&&S^{10}\smsh S^{10}\ar[u]\\
}$$
where the bottom square comes from the decomposition of
$\d_2$ through the Hopf invariant \cite{FGT}.
The space
$\overline{X}/\overline{A}= (S^{10}\cup D^{14}\cup D^{18})\cup D^{21}$
is of cone length less than 2 since $S^{10}\cup D^{14}\cup D^{18}$ is a suspension. 
Denote by
$\varphi
\colon S^{20}\rightarrow S^{10}\cup D^{14}\cup D^{18}$ the attaching map of the 20-cell. The
composite of $\varphi$  with the projection $S^{10}\cup D^{14}\cup D^{18}\rightarrow
S^{14}\cup D^{18}$ is stably trivial, therefore trivial. As a consequence, the map $\varphi$ lifts
to
$S^{10}$. The morphism of suspension
$E\colon \pi_{19}(S^9)\rightarrow \pi_{20}(S^{10})$
is surjective \cite[Theorem\thinspace 7.3]{Toda}, 
and therefore the Hopf invariant
$\overline{H}(\varphi)$ is zero.

%(In fact, this proves that $X/A$ is a suspension!)

\smallskip

\emph{Proof of {\em (b)}}
From the product
structure of $H^*(\sp(3);\ZZ)$ we deduce that 
the map $\alpha:S^{11}\wdg S^{10}
\maprt{} S^{10}\wdg S^{10}$ is given by
$\alpha = ([ b\eta, c\eta], [\pm 1,\pm 1] ) $
for some constants $b,c\in \ZZ/2$.
We prove now that $b=c=1$.

Consider the following diagram obtained by collapsing subcomplexes,
cellular approximations and 
the decomposition of $\d_2$ through Hopf-Ganea 
invariants:
$$\xymatrix{
S^3\cup D^7\cup D^{11}\ar[r]^-{\d_2}\ar[d]&
C\smsh C\ar[rrr]^{1\smsh \d_2}\ar[d]\ar[dr]&&&
C\smsh C^{\smsh 2}\\
S^{11}\ar[r]\ar[rd]_-{[b\eta,c\eta]}&
(C\smsh C)/S^6\ar[r]&C\smsh S^7\ar[rr]^{1\smsh \eta}&&
C\smsh S^6\ar@{^{(}->}[u]\\
&S^{10}\smsh S^{10}\ar@{^{(}->}[u]\ar[ur]^{\gamma}
\ar[rrr]_{\overline{\gamma}}&&&S^9.\ar@{^{(}->}[u]
}$$
The map $\gamma$  induces the homomorphism
$[1,0]$ on $H^{10}$ and hence also on the homotopy groups
$\pi_{10}$, which means that $\overline{\gamma}=[\eta,0]$, and the
composite $S^{11}\rightarrow S^9$ is $b\eta^2$.
Since this composite is known to be nonzero \cite{FGT},
it follows that $b=1$.

The same proof, but now using
$C\smsh C\rightarrow S^7\smsh C$,
reveals that $c=1$.

\smallskip

\emph{Proof of {\em (c)}}
The map $\beta: S^{10}\wdg S^{11}\maprt{} S^9$ must have
the form $\beta = (s\eta, t\eta^2)$ with $s,t\in \ZZ/2$.  It is defined by the
diagram (with auxiliary space $Z$ to be discussed below):
$$
\xymatrix{Z\ar[r]^i \ar[ddd] &
\overline X\ar[rr]^{\overline\Delta_3} \ar[d]
   && \overline X^{\smsh 3}\ar@{=}[d]\\
&\overline X/\overline X_7  \ar[rr] && \overline X^{\smsh 3}\\
&(S^{10}\wdg S^{11} )\cup D^{14}
    \ar[rr] \ar@{^{(}->}[u]_{17-\mathrm{skeleton}}
  && C^{\smsh 3} \ar@{^{(}->}[u]_{15-\mathrm{skeleton}}\\
S^{k} \ar[r]^(.4)j
  & S^{10}\wdg S^{11}\ar[rr]^(.58)\beta \ar@{^{(}->}[u]
  && S^9
\ar@{^{(}->}[u]_{12-\mathrm{skeleton}}.\\
}
$$
First we show that the first coordinate of $\beta$ is $\eta$.
For this, we choose $Z = \sp(2)\sseq X$ (and $k=10$).  According to
Schweitzer \cite{Schw}, $\wcat(\sp(2)) = 3$.  Since
the inclusion $\sp(2)^{\smsh 3}\maprt{} X^{\smsh 3}$ is a $12$-equivalence,
we conclude that the bottom composite is essential,
and this identifies the first component.

We study the second coordinate by taking $k=11$ and
$Z = S^{3}\cup D^{7}\cup D^{11}\sseq \overline X$.
By cellular approximation there is a unique map
$Z\maprt{} S^9$ so that the composite $Z\maprt{}S^9\inclds X^{\smsh 3}$
is $\overline\Delta_3 |_Z$.  Since $[\s C, S^9 ]=[C,S^9] =*$, the
quotient $Z\maprt{}S^{11}$ induces
an isomorphism $[S^{11},S^9]\maprt{\cong} [Z,S^9]$.
From the fact that $\wcat(Z)=3$ (cf. Proposition \ref{prop:cone}), we deduce
that the
composite $S^{11} \maprt{} S^{9}$ is $\eta^2$.
\hfill$\Box$\par

With exactly the same argument than in Part (a) of Lemma~\ref{lem:pab},
one can prove:

\blem\label{lem:pabb}
The composite map
$X\rightarrow X \smsh X\rightarrow X/A\smsh X/A$
is trivial.
\elem

\section{Proof of Proposition \ref{prop:C}}

We may take the homotopy equivalence $T$ to be 
induced by an equivalence $\twdl T$ which makes the diagram 
$$
\xymatrix{
S^{19}\wdg S^{19}\wdg S^{19} \ar[rr]^{\twdl T}
     \ar[d]_{(\nu,2\nu,2\nu)}
   && S^{19}\wdg S^{19}\wdg S^{19} \ar[d]^{(\nu,0,0)}
\\
S^{16} \ar@{=}[rr] \ar[d] && S^{16}  \ar[d] 
\\
{(C\smsh C\smsh C_{10})}_{22} \ar[rr]^{T}_{\simeq} && 
   C_{16}\wdg S^{20}\wdg S^{20}  }
$$
commute.  There are many maps $\twdl T$ which suffice
for this purpose; for simplicity, we choose to 
work with $\twdl T = ([1,0,0], [-6,1,-3],[2,0,1] )$.
 
We now consider the diagram
$$\xymatrix{
S^{21} \ar[rr]^(.4)q \ar[rrd] 
   && C\smsh C\smsh C_{10} \ar[rr]^{T}\ar[d]^{u} 
   && C_{16}\wdg S^{20}\wdg S^{20} \ar[d] 
      \ar[r]& C^{\smsh 5}\ar[d]^v
\\
   && S^{20}\wdg S^{20}\wdg S^{20} \ar[rr]^{\s\twdl T}
   && S^{20}\wdg S^{20}\wdg S^{20} \ar[r] & C^{\smsh 5}/S^{15} }
$$
in which $\s\twdl T$ has the same explicit description
as $\twdl T$.

\smallskip

\emph{Proof of {\em (a)}}  From the explicit
description of $\s\twdl T$ we see that it is
sufficient to show $u\of q = [\lambda, 0, 0] :
S^{21}\maprt{} S^{20}\wdg S^{20}\wdg S^{20}$.  
Recall from Section 3 that $A$ denotes 
the subcomplex $S^3 \cup \c (C_6)\sseq X$, which has
weak category~$2$.  The diagram
$$\xymatrix{
S^{21} \ar[rr] \ar[d] 
   && X/X_3 \smsh X/A  \ar[rr]\ar[d] 
   && (X\smsh X)\smsh X/A  \ar[d] 
\\
S^{21} \ar[rr]^-{*}\ar[rrd]
   && X/A \smsh X/A  \ar[rr]
   && (X\smsh X/ S^{3}\smsh S^{3})\smsh X/A 
\\
&&  S^{20}\ar[rr]^{[1,1]}\ar@{^{(}->}[u]&& S^{20}\wdg S^{20}\ar@{^{(}->}[u] }
$$
shows that the composite 
$$
S^{21} \maprt{} C\smsh C\smsh C_{10}
\maprt{} S^{20}\wdg S^{20}\maprt{}(X\smsh X)/ (S^{3}\smsh S^{3})\smsh X/A
$$
factors through $X\maprt{} X/A\smsh X/A$
which was shown to be trivial in Lemma~\ref{lem:pabb}.  
Since $S^{20}\wdg S^{20}\maprt{}(X\smsh X/ S^{3}\smsh S^{3})\smsh X/A$ is a 
$22$-equivalence, this proves (a).

\smallskip

\noindent{{\sc Proof of} (b)}\quad   The map we are interested
in is the composite along the bottom of the diagram
$$\xymatrix{
S^{21} \ar[rr] \ar[d] 
   && X/X_3 \smsh X/A  \ar[rr]\ar[d] 
   && (X\smsh X)\smsh X/A  \ar[d] 
\\
S^{21} \ar[rr]^-{*}\ar[rrd]_{\simeq}
   && X/X_3 \smsh X/X_{11}  \ar[rr]
   && (X\smsh X)\smsh X/X_{11} 
\\
&&  S^{7}\smsh S^{14}\ar[rr]^{\eta\smsh 1}\ar@{^{(}->}[u]
   && S^{6}\smsh S^{14}\ar@{^{(}->}[u] }
$$
%$$
%\xymatrix{
%X\ar[rr]\ar@{=}[d] && X\smsh X\ar[d]  \ar[rr]^{\d_2\smsh 1} 
%    && (X\smsh X)\smsh X \ar@{=}[d]
%\\
%X\ar[rr]\ar[d] && (X/X_3) \smsh X\ar[d]  \ar[rr]  
%    && (X\smsh X)\smsh X  
%\\
%S^{21}\ar[rr] \ar@{=}[d] && S^7\smsh S^{14}\ar[rr]\ar@{=}[d]  
%    && (X\smsh X)\smsh S^{14} \ar@{^{(}->}[u]
%\\
%S^{21}\ar[rr]^{\simeq}  && S^7\smsh S^{14}\ar[rr] 
%    && (S^3\smsh S^3)\smsh S^{14} \ar@{^{(}->}[u]
%}
%$$
The cup product structure of $\sp(3)$ shows that the
map $S^{21}\maprt{} S^{7}\smsh S^{14}$ is a homotopy
equivalence; and the map $S^{7}\smsh S^{14}\maprt{}
S^{6}\smsh S^{14}$ is evidently the $14$-fold
suspension of the map $S^{7}\maprt{}
S^{6}$ induced by $\d_2$.
This latter map
was shown in \cite{Schw} to be $\eta$.

Finally, we determine the composite $C_{16}\maprt{}
C^{\smsh 5}\maprt{} C^{\smsh 5}/S^{15}$.  
The $22$-skeleton of $C^{\smsh 5}/S^{15}$ is 
$(\bigvee_{i=1}^5 S^{19})\wdg (\bigvee S^{22})$.
Therefore this composite factors through 
$$
S^{20}\maprt{w}\bigvee_{i=1}^5 S^{19}\inclds
C^{\smsh 5}\maprt{} C^{\smsh 5}/S^{15},
$$ 
where $w=[\epsilon_i \eta ]$ is either $0$ or $\eta$ in each coordinate.  
Consider the commutative diagram
$$
\xymatrix{
S^{21} \ar[rr]^{\lambda} \ar[rrd]_{\eta} \ar@(ur,ul)[rrrr]^{h} %\ar@/^1pc/[rrrr]^{h}
   && C_{16} \ar[d] \ar[rr] 
   && C^{\smsh 5}\ar[d]^v
\\
&& S^{20} \ar[rr] && C^{\smsh 5}/S^{15}, }
$$
with $v\of h\simeq *$ by Proposition \ref{prop:B}.
The composite can also be identified with $[\epsilon_i \eta^2 ]$,
which shows that each $\epsilon_i = 0$ and therefore the map 
$S^{20}\maprt{} C^{\smsh 5}/S^{15}$ is trivial.
\hfill$\Box$\par

%%%%%%%%%%%%%%%%%%%%   End of main body of article
%
%                             References
%

\small

Centro de Matem\'atica (CMAT)\\
Universidade do Minho (Gualtar)\\
4710 Braga, Portugal\\
\texttt{lfernandez@math.uminho.pt}\\

Departamento de Xeometr\'{\i}a e Topolox\'{\i}a\\
Universidade de Santiago de Compostela\\
15706 Santiago de Compostela\\
Espa\~na\\
\texttt{agtato@zmat.usc.es}\\

Department of Mathematics,
Dartmouth College\\
Hanover, NH 03755,
U.S.A.\\
\texttt{Jeffrey.A.Strom@Dartmouth.EDU}\\

D\'epartement de Math\'ematiques\\
UMR 8524\\
Universit\'e de Lille~1\\
59655 Villeneuve d'Ascq Cedex,
France\\
\texttt{Daniel.Tanre@agat.univ-lille1.fr}\\

%
%
%
%\vfill\eject
%%
%    Type any appendix material (to go after the references) here
%
\end{document}